\documentclass[a4paper,reqno,11pt]{amsart}

\usepackage[latin1]{inputenc}
\usepackage[T1]{fontenc}

\usepackage{amssymb,amsmath}
\usepackage{setspace}
\usepackage[normalem]{ulem}
\usepackage{color}
\usepackage{geometry}
\usepackage[english]{babel}

\newtheorem{proposition}{Proposition}
\newtheorem{lemma}{Lemma}

\newcommand{\highlight}[1]{#1}
\newcommand{\highlightt}[1]{#1}

\newcommand{\M}{T}
\def\eps{\varepsilon}
\def\WVG{WMG}

\begin{document}
\title{
The nucleolus of large majority games}\thanks{We thank the Maria Montero, Roberto Serrano and an anonymous referee for helpful discussion and constructive comments on an earlier version.}

\author{Sascha Kurz}
\address{Sascha Kurz, Department of Mathematics, University of Bayreuth, Germany. Tel: +49 921 557353,
Fax: +49 921 557352. E-mail: sascha.kurz@uni-bayreuth.de}

\author{Stefan Napel}
\address{Stefan Napel, Department of Economics, University of Bayreuth, Germany}

\author{Andreas Nohn}
\address{Andreas Nohn, Public Choice Research Centre, University of Turku, Finland}


\begin{abstract}
\small Members of a shareholder meeting or legislative committee have greater or smaller voting
power than meets the eye if the nucleolus of the induced  majority game differs from the voting
weight distribution. We establish a new sufficient condition for the weight and power distributions
to be equal; and we characterize the limit behavior of the nucleolus in case all relative weights become small.


\medskip


\medskip

\end{abstract}

\maketitle

\setstretch{1.2}

\section{Introduction}
\label{sec:intro}
\noindent

\noindent Among all individually rational and efficient payoff vectors in a
game $v$ with transferable utility, the \emph{nucleolus} selects
a particularly stable one. It quantifies each coalition's dissatisfaction with
a proposed vector $x$ as the gap between the coalition's worth $v(S)$ and the
surplus share $\sum_{i\in S} x_i$ that is allocated to  members of $S\subseteq N$;
then it selects the allocation $x^*$ which involves lexicographically minimal dissatisfaction.
In contrast to other prominent point solutions in cooperative game theory, such as the Shapley
value, $x^*$ is guaranteed to lie in the \emph{core} of game $(N,v)$ whenever that is non-empty.

Even before the final version of Schmeidler's article\nocite{Schmeidler:1969} which established the
definition, existence, uniqueness, and continuity of the nucleolus was published in 1969,
Peleg, see \cite{Peleg:1968}, had applied it to \emph{weighted majority games (\WVG)}.
In these games the worth of a coalition $S$ of players is either 1 or 0, i.e., $S$ is either winning or
losing, and there exists a non-negative quota-and-weights \emph{representation} $[q; w_1, \ldots, w_n]$
such that $v(S)=1$ iff $\sum_{i\in S}w_i\ge q$.
The weight vectors that constitute a representation of a given {\WVG} $v$ for some quota $q$ form a
non-singleton convex set $R(v)$.

Peleg highlighted a property of constant-sum {\WVG}s with a \emph{homogeneous} representation, i.e.,
one where total weight of any minimal winning coalition equals $q$: the nucleolus $x^*$ of such
a {\WVG} $v$ is contained in $R(v)$, i.e., it is also a representation.\footnote{
A {\WVG}
$(N,v)$ is called \emph{constant-sum} if for any $S\subseteq N$ either $v(S)=1$ or $v(N\smallsetminus S)=1$.
$S\subseteq N$ is a \emph{minimal winning coalition (MWC)} if $v(S)=1$ and $v(T)=0$ for any $T\subset S$.}
Despite this early start, the relation between voting weights and the nucleolus of weighted majority games~--
constant-sum or not, homogeneous or inhomogeneous~-- has to the best of our knowledge not been studied systematically
so far. This paper is a first attempt to fill this gap.

Discrepancies between weights and the nucleolus matter because the nucleolus is an important indicator of
influence in collective decision bodies. It emerges as an equilibrium price vector in models that evaluate
voters' attractiveness to competing lobbying groups, see \cite{Young:1978,Shubik/Young:1978};
more recent theoretical work by Montero (\cite{Montero:2005,Montero:2006}) has established it as a focal
equilibrium prediction for strategic bargaining games with a majority rule.\footnote{Corresponding experimental lab evidence is
mixed; see \cite{Montero/Sefton/Zhang:2008}. Non-cooperative foundations of the nucleolus for other than majority games have been given, e.g., by \cite{Potters/Tijs:1992} and \cite{Serrano:1993,Serrano:1995}.}
So large differences between a voter~$i$'s weight $w_i$ and nucleolus $x^*_i$ can mean that the real power
distribution in a decision body such as a shareholder meeting is hidden from the casual observer.
This intransparency can be particularly problematic for political decision bodies, where voting weight
arrangements affect the institution's legitimacy.\footnote{See \cite{LeBreton/Montero/Zaporozhets:2012} for
nucleolus-based power analysis of the European
Union's Council; an early-day weight arrangement meant that Luxembourg had a relative voting weight of 1/17 but
zero voting power.~-- In general, the power-to-weight ratio can differ arbitrarily from 1. For instance, the
nucleolus of the {\WVG} with representation $[0.5; (1-\eps)/2, (1-\eps)/2, \eps]$ is
$x^*=(1/3, 1/3, 1/3)$ for any $\eps\in (0; 0.5)$.}

This paper investigates absolute and relative differences between players' relative voting weights as defined
by vote shares in an assembly, electoral college, etc.\ and the nucleolus of the implied {\WVG}. We determine an
upper bound on their $\Vert\cdot\Vert_1$-distance which depends only on
quota and maximum weight in a given representation in Lemma~\ref{lemma_distance}.
The lemma allows to conclude that if the relative weight of every individual voter in player set
$\{1, \ldots, n\}$ tends to zero, then the ratio $x_i^*/x_j^*$ of two nucleolus components converges to
$w_i/w_j$ for all regular voters $i$ and $j$ as $n\to \infty$ (Prop.~\ref{thm_distance_convergence}).
This complements analogous limit results in the literature on the Shapley value, the Banzhaf value and
voter pivotality on intervals, see \cite{Kurz/Maaser/Napel:2013,Lindner/Machover:2004,Neyman:1982}, as
well as for stationary equilibrium payoffs in legislative bargaining games \`{a} la Baron-Ferejohn, see 
\cite{Snyder/Ting/Ansolabehere:2005}. We also establish a new sufficient condition for the
nucleolus to coincide with given relative weights (Prop.~\ref{main_thm}).
It implies that a finite number of replications brings about full coincidence for any given {\WVG}.

\section{Nucleolus}
\label{sec:nucleolus}
\noindent
Consider a {\WVG} $(N,v)$ with representation $[q; w_1, \ldots, w_n]$. Using  notation
$x(S)=\sum_{i\in S}x_i$, a vector $x\in\mathbb{R}^n$ with $x_i\ge v(\{i\})$ and $x(N)=v(N)$ is called an
\emph{imputation}. For any coalition $S\subseteq N$ and imputation $x$, call $e(S,x)=v(S)-x(S)$ the
\emph{excess} of $S$ at $x$. It can be interpreted as quantifying the coalition's dissatisfaction and
potential opposition to an agreement on allocation $x$.
For any fixed $x$ let $S_1, \ldots, S_{2^n}$ be an ordering of all coalitions such that the excesses at
$x$ are weakly decreasing, and denote these ordered excesses by $E(x)=\big(e(S_k,x)\big)_{k=1, \ldots, 2^n}$.
Imputation $x$ is \emph{lexicographically less} than imputation $y$ if $E_k(x)<E_k(y)$ for the smallest
component $k$ with $E_k(x)\neq E_k(y)$.
The \emph{nucleolus} of $(N,v)$ is then uniquely defined as the lexicographically minimal
imputation.\footnote{Schmeidler's (\cite{Schmeidler:1969}) original definition did not restrict the
considered vectors to be imputations but is usually specialized this way. The set of imputations that minimize
just the largest excess, $E_1(x)$, is called the \emph{nucleus} of $(N,v)$ by Montero \cite{Montero:2006}. Our results
are stated for the nucleolus but apply to every element of the nucleus: 
both coincide under the conditions of
Prop.~\ref{main_thm}; Lemma~\ref{lemma_distance} and Prop.~\ref{thm_distance_convergence} generalize straightforwardly.}

As an example, consider $(N,v)$ with representation $[q;w]=[8;6,4,3,2]$. The nucleolus can be computed as
$x^*=(2/5, 1/5, 1/5, 1/5)$ by solving a sequence of linear programs~-- or by appealing to the sufficient condition
of \cite{Peleg:1968} after noting that the game is constant-sum and permits a homogeneous representation $[q';w']=[3;2,1,1,1]$.
Denoting the \emph{normalization} of weight vector $w$ by $\bar w$, i.e., $\bar w=w/\sum w_i$, the respective total
differences between relative weights and the nucleolus are $\Vert \bar w-x^*\Vert_1=2/15$ for the first and
$\Vert \bar w'-x^*\Vert_1=0$ for the second representation (with $\Vert x\Vert_1=\sum |x_i|$).

\section{Results}
\label{sec:distance}

\noindent Saying that representation $[q;w]$ is \emph{normalized} if $w=\bar w$, we have:\footnote{All proofs are provided
in the Mathematical Appendix.}
\begin{lemma}
\label{lemma_distance}
\doublespacing
Consider a normalized representation $[q;w]$ with $0<q<1$ and $w_1\geq\dots\geq w_n\geq 0$ and let $x^*$ be the
nucleolus of this {\WVG}. Then
\begin{equation}\label{eq:lemma}
\Vert x^*-w\Vert_1 \leq \frac{2w_1}{\min\{q,1-q\}}.
\end{equation}
\end{lemma}

If we consider a sequence $\{(\{1, \ldots, n\}, v^{(n)})\}_{n\in\mathbb{N}}$ of $n$-player {\WVG}s with
representations $[q^{(n)};w^{(n)}]$ such that the normalized quota $\bar q^{(n)}$ is bounded away from 0 and 1 (or,
more generally, 0 and 1 are no cluster points of $\{\bar q^{(n)}\}_{n\in\mathbb{N}}$),
and each player~$i$'s normalized weight $\bar w_i^{(n)}$ vanishes as $n\to \infty$ then Lemma~\ref{lemma_distance} implies
\begin{equation}
\lim_{n\to\infty} \Vert x^{*(n)}-\bar w^{(n)} \Vert_1\to 0.
\end{equation}

Convergence to zero of the total difference between nucleolus components $x^{*(n)}_i$ and relative voting weights
$\bar w^{(n)}_i$ does not yet guarantee that the nucleolus is asymptotically proportional to the weight vector, i.e.,
that each ratio $x^{*(n)}_i/x^{*(n)}_j$ converges to $w_i/w_j$.
This can be seen, e.g., by considering
\begin{equation}\label{eq:pathological_example}
    \big[q^{(n)};w^{(n)}\big]=\Bigg[\frac{2n-1}{2};1, \underbrace{2, \ldots, 2}_{n-1}\Bigg].
\end{equation}
The nucleolus either equals $\big(0,\frac{1}{n-1},\dots,\frac{1}{n-1}\big)$ or $\big(\frac{1}{n},\dots,\frac{1}{n}\big)$
depending on whether $n$ is even or odd; ratio $x^{*(n)}_1/x^{*(n)}_2\neq \frac{1}{2}$ alternates between 0 and 1.

But such pathologies are ruled out for players $i$ and $j$ whose weights are ``non-singular'' in the weight sequence
$\big\{w^{(n)}\big\}_{n\in\mathbb{N}}$. Specifically, denote the total number of players~$i\in \{1,\ldots, n\}$ with
an identical weight of $w^{(n)}_i=\omega$ by $m_\omega(n)$.
We say that a player~$j$ with weight $w_j$ is \emph{regular} if $m_{w_j}(n)\cdot \bar w^{(n)}_j$ is bounded away
from 0 by some constant $\eps>0$.
Lemma~\ref{lemma_distance} then implies:\footnote{We assume $w^{(n)}_j=w_j$ in our exposition. Adaptations to cases
where $q^{(n)}$ \emph{and} $w^{(n)}_j$
vary in $n$ are straightforward. The essential regularity requirement is that a voter type's aggregate relative weight
does not vanish.}
\begin{proposition} \label{thm_distance_convergence}
Consider a sequence $\big\{\big[q^{(n)};(w_1, \ldots, w_n)\big]\big\}_{n\in\mathbb{N}}$ with corresponding normalized
quotas that exclude 0 and 1 as cluster points and with normalized weights satisfying $\bar w^{(n)}_k\downarrow 0$ for
every $k\in \mathbb{N}$ as $n\to \infty$.
Then the nucleolus $x^{*(n)}$ of the {\WVG} represented by $\big[q^{(n)};(w_1, \ldots, w_n)\big]$ satisfies
\begin{equation}
\lim_{n\to \infty} \frac{x^{*(n)}_i}{x^{*(n)}_j}=\frac{w_i}{w_j}
\end{equation}
for any regular players $i$ and $j$.
\end{proposition}

For a considerable class of games, asymptotic equality of nucleolus and normalized weights can be strengthened to
actual equality.\footnote{Non-null players have a positive nucleolus value in this case~-- in contrast to {\WVG}s
in general. So we implicitly establish a sufficient condition for $w_i>0 \Rightarrow x_i^*>0$.} Namely, for a fixed
$n$-player {\WVG} with representation $[q; w_1, \ldots, w_n]$ let $m_{\omega}$ denote the number of players
that have weight $\omega$; so
\begin{equation}
m^\circ=\min_{i\in\{1,\ldots, n\}} m_{w_i}\ge 1
\end{equation}
is the number of occurrences of the rarest weight in vector $w=(w_1, \ldots, w_n)$.

\begin{proposition}\label{main_thm}
Consider a {\WVG} representation $[q; w]$ with integer weights $w_1\ge \ldots \ge w_n\ge 0$ and normalization
$[\bar q;\bar w]$. Denoting the number of numerically distinct values in $w$ by $1\le t \le n$, the nucleolus $x^*$ of
this {\WVG} satisfies
\begin{equation}\label{eq:main_thm}
x^*=\bar w \text{\ \ if \ }
\min\{\bar q, 1-\bar q\}\cdot m^\circ > 2t {w_1}^2.
\end{equation}
\end{proposition}

The proposition refers to \emph{integer} weights. Even though it is not difficult to obtain an integer representation
from any given $[q;w]$ with non-integer values, this is an important restriction. In particular, it is not possible to
rescale a given weight vector $w$ so as to make the maximal weight $w_1$ arbitrarily small.\footnote{Note also that
inequality~(\ref{eq:main_thm}) must be violated if two interchangeable players of $(N,v)$ have different weights because
$x^*=\bar w$ would then contradict the symmetry property of the nucleolus. So as a subtle implication of the integer
requirement, weight changes which would destroy a given symmetric or `type-preserving' representation and satisfy
(\ref{eq:main_thm}) are impossible.
Another way to look at this is to begin with a {\WVG}'s representation where $w_i\neq w_j$ for interchangeable players~$i$ and $j$ and then to replicate all players and weights:
after enough replications the two players (types) $i$ and $j$ must lose their interchangeability.}

The right-hand side of the inequality in condition~(\ref{eq:main_thm}) is smaller, the smaller the number of different
weights in the representation, and the smaller the involved integers (particularly $w_1$). Similarly, the left-hand
side is larger, the greater the number of occurrences of the rarest weight.
It follows that condition~(\ref{eq:main_thm}) is most easily met when null players (where $x_i^*=0$ is known) are
removed from the {\WVG} in question and a \emph{minimal integer representation} is considered.\footnote{Uniqueness
and other properties of minimal integer representations of {\WVG} are investigated in \cite{Kurz:2012:representation}.}
This is automatically also a homogeneous representation if any exists.

Our sufficient condition for $x^*=\bar w$ is, however, independent of the known homogeneity-based one. The normalization of
weights in $[3;2,1,1,1]$ must, according to Peleg \cite{Peleg:1968}, coincide with the {\WVG}'s nucleolus because the
game is constant-sum; but our condition~(\ref{eq:main_thm}) is violated.
In contrast, the representation $[q;w]=[1500; 4,\dots,4, 3,\dots,3, 2,\dots,2]$ of a 900-player {\WVG} where each
of the $t=3$ weight types occurs $m^\circ=300$ times satisfies our condition.
Hence $x^*=\bar w$. Since the game is inhomogeneous,\footnote{Coalitions with (a) 300, 100, and 0, (b) 300, 0,
and 150, or (c) 300, 1 and 149 players of weights 4, 3, and 2 are minimal winning, and cannot be made to have identical
aggregate weights in any representation $[q'; w']$.} Peleg's finding does not apply.

The left-hand side in condition~(\ref{eq:main_thm}) equals at most half the number of occurrences of the rarest weight,
$m^\circ$, and the right-hand side is bounded below by 2.
This, first, implies that the condition cannot be met by {\WVG}s where only one instance of some weight type is
involved. This limits Prop.~\ref{main_thm}'s applicability for small-scale games such as $[3;2,1,1,1]$.
But, second, it means that if we consider \emph{$\rho$-replicas} of any given $n$-player {\WVG} with integer
representation $[q; w]$, i.e., a {\WVG} with quota $\rho q$ and $\rho$ instances of any of the $n$ voters in
$[q; w]$, then one can compute an explicit number $\tilde\rho$ from (\ref{eq:main_thm}) such that the nucleolus of the
resulting $\rho n$-player {\WVG} must coincide with the corresponding normalized weight vector for every
$\rho\ge \tilde \rho$.\footnote{For simple majority games which involve equal numbers of voters with weight $\omega=4$,
3, and 2, condition~(\ref{eq:main_thm}) calls for $m^\circ > 192$. But $x^*=\bar w$ already holds after one replication
of $[5; 4,3,2]$, i.e., $\rho\ge 2$. So tighter bounds might be obtained by different techniques than ours. However,
surprising sensitivity of $x^*$ to the game at hand cautions against too high expectations. For instance, $x^*=\bar w$
if $w_1=\ldots=w_5=4$ and either $w_6=\ldots=w_{11}=1$ or $w_6=\ldots=w_{13}=1$ with $\bar q=58\%$; in contrast, $x^*=(1/5, \ldots, 1/5, 0, \ldots, 0)$ if $w_6=\ldots=w_{12}=1$. We thank Maria Montero for suggesting this example.}
This observation echoes the coincidence result obtained in \cite{Snyder/Ting/Ansolabehere:2005} for {\WVG} replicas
under Baron-Ferejohn bargaining:\footnote{We thank an anonymous referee for pointing out to us that Snyder et al.'s Prop.~2 is in fact a corollary to our Prop.~\ref{main_thm}, the uniqueness of SSPE payoffs recently established in \cite{EraslanMcLennan:2013:uniqueness}, and Montero's (2006) Prop.~7.}
at least in sufficiently large majority games, voting weight and power are the same.


\section*{Mathematical Appendix}

\noindent
\textbf{Proof of Lemma~\ref{lemma_distance}.}
\highlightt{Let $P$ be the set of players $i$ such that $\{i\}$ is a winning coalition of $[q;w]$. If $|P|>1$, then 
the nucleolus is not defined. If $|P|=1$, then $w_1\ge q$, so that the upper bound is trivially satisfied due to 
$\Vert x-w\Vert_1\le 2$ for any $x\in\mathbb{R}_{\ge 0}^n$ with $\Vert x\Vert_1=1$. In the remaining cases we have $x^*(S)\ge q$ 
for any winning coalition, since chosing $x=w$ yields a maximum excess $\max_{S\subseteq N} E(S,x)$ of at most $1-q$.}

Define $w(S)=\sum_{i\in S}w_i$ and $x^*(S)=\sum_{i\in S}x^*_i$. Let $S^+=\{i\in N \mid x^*_i>w_i\}$ and
$S^-=\{i\in N \mid x^*_i\leq w_i\}$.
\highlight{We have $w(S^+)<1$ since $w(S^+)<x(S^+)\le x(N)=1$, so that $w(S^-)>0$.}
\highlight{Define} $0\leq \delta\leq 1$ by $x^*(S^-)=(1-\delta)w(S^-)$. 
We have 
  \begin{equation}
    \highlight{x^*(S^+)=1-x^*(S^-)=w(S^+)+w(S^-)-(1-\delta)w(S^-)=w(S^+)+\delta w(S^-)}
  \end{equation}   
  and
  \begin{equation}
    \label{eq:key_lemma_observation2}
    \highlight{\Vert w\!-\!x^*\Vert_1\!=\!\sum\limits_{i\in S^+} \!\!\!\left(x^*_i\!-\!w_i\right)\!+\sum\limits_{i\in S^-} \!\!\!\left(w_i\!-\!x^*_i\right)\!
    =x^*(S^+)\!-\!w(S^+)\!+\!w(S^-)\!-\!x^*(S^-)=2\delta w(S^-).}
  \end{equation}
Let $T$ \highlight{be} 
\highlight{generated} by starting with $S=\varnothing$ and successively adding a remaining player~$i$ with
minimal $x^*_i/\highlightt{w_i}$ until $w(T)\ge q$.

In case $w(S^-)\geq q$ we then have $x^*(T)/w(T)\leq x^*(S^-)/w(S^-) =1-\delta$. Multiplying by $w(T)$, using
$q\leq w(T)\leq q+w_1$ and finally $\delta \le 1$ yields
  \begin{equation}
    x^*(T)\leq (1-\delta)w(T)\leq (1-\delta)(q+w_1)\leq q(1-\delta)+w_1.
  \end{equation}
This and \highlight{$x^*(T)\geq q$} 
yield $\delta\leq w_1/q$.
Applying this and $w(S^-)\highlight{\leq} 1$ in equation~(\ref{eq:key_lemma_observation2}) gives $\Vert x^*-w\Vert_1\leq \frac{2w_1}{q}$.

In case $w(S^-)<q$ \highlight{we have $w(S^+)>0$. Note} that moving from $S^-$ to $T$ involves the weight addition $w(T)-w(S^-)$ which comes with a
nucleolus per weight unit of at most $x^*(S^+)/w(S^+)$. So
\begin{eqnarray}
  x^*(T)& = & x^*(S^-)+x^*(T\backslash S^-)\notag\\
  &\le& (1-\delta)w(S^-) + \frac{x^*(S^+)}{w(S^+)}\cdot \bigl(w(T)-w(S^-)\bigr)\notag\\
  &\le& (1-\delta)w(S^-) + \frac{x^*(S^+)}{w(S^+)}\cdot \bigl(q-w(S^-)+w_1\bigr)
\end{eqnarray}
where the last inequality uses $w(T)\leq q+w_1$. Rearranging with $x^*(S^+)=w(S^+)+\delta w(S^-)$ and  $w(S^-)=1-w(S^+)$ yields
\begin{eqnarray}
    x^*(T) &\leq& q+\frac{w(S^+)+\delta w(S^-)}{w(S^+)}\cdot w_1-\frac{(1-q)\delta w(S^-)}{w(S^+)}. \label{ie_pl112}
\end{eqnarray}
Since $\delta\le 1$ the right hand side of (\ref{ie_pl112}) is at most $q +\bigl(w_1 - (1-q)\delta w(S^-)\bigr)/w(S^+)$.
So $q\leq x^*(T)$ implies $(1-q)\delta w(S^-)\leq w_1$.  Hence $\Vert x-w\Vert_1 \leq \frac{2w_1}{1-q}$.
\hfill{$\square$}

\bigskip

\noindent
\textbf{Proof of Proposition~\ref{thm_distance_convergence}.}
If $x_i^{*(n)}/\,\bar{w}_i^{(n)}\ge 1+\delta$ or $x_i^{*(n)}/\,\bar{w}_i^{(n)}\le 1-\delta$ then
 $\Vert x^{*(n)}-\bar{w}^{(n)}\Vert_1\ge \delta \cdot m_{w_i}(n)\cdot \bar{w}_i^{(n)}\ge \delta\varepsilon$ for
 some $\varepsilon>0$ if $i$ is regular.
But $\lim\limits_{n\rightarrow\infty} ||x^{*(n)}-\bar w^{(n)}||_1=0$ (Lemma~\ref{lemma_distance}). So $\lim\limits_{n\rightarrow\infty} {x^{*(n)}_i}/\,{\bar w^{(n)}_i}=1$ and hence
\begin{equation}
   1=\lim_{n\to\infty} \frac{x_i^{*(n)}}{\bar{w}_i^{(n)}}\cdot \frac{\bar{w}_j^{(n)}}{x_j^{*(n)}}=
   \lim_{n\to\infty} \frac{x_i^{*(n)}}{x_j^{*(n)}}\cdot \frac{w_j}{w_i} \text{\ \ \ if $i$ and $j$ are regular.}
\end{equation}

\hfill{$\square$}

\noindent
\textbf{Proof of Proposition~\ref{main_thm}.}
It suffices to prove the result in case $w_n>0$ because $w_i=0$ directly implies $x^*_i=0$. We may also assume $0<\bar q<1$.
For each $k\in \{1, \ldots, t\}$ let $\omega_k$ denote the normalized weight of a voter~$i$ with type~$k$ (i.e., $\bar w_i=\omega_k$) and, with slight abuse of notation, let $x^*_k$ be this voter/type's nucleolus.
Define $r_k=x^*_{k}/\omega_{k}$ and w.l.o.g.\ assume $r_1\ge \ldots \ge r_t$.
Let $a$ denote the largest index such that $r_1=r_a$ and $b$ be the smallest such that $r_b=r_t$.
The claim is true if $a\ge b$.
So we suppose $a<b$ and establish a contradiction by showing that we can construct an imputation $x^{**}$ with maximum excess $E_1(x^{**})$ smaller than $E_1(x^{*})$.

Writing $\varepsilon=\frac{1}{2}\cdot \min\{\bar q,1-\bar q\}$ and $n_k=m_{\omega_k}$, the premise and $t, w_1\ge 1$ imply
\begin{equation}\label{ie_intermed_step}
   w_1\le t {w_1}^2 <\eps m^\circ \leq \eps n_k
\end{equation}
for each $k\in \{1, \ldots, t\}$. Considering $\omega_k$-weighted sums of (\ref{ie_intermed_step}) we obtain
\begin{equation*}
  \hspace{-0.55cm}\text{(I)} \quad  \sum_{k\le a} w_1 \omega_k < \varepsilon\sum_{k\le a} n_k \omega_k \quad \text{\ and \ \quad \ (II)}\quad \sum_{k\ge b} w_1 \omega_k < \varepsilon\sum_{k\ge b} n_k \omega_k. \notag
\end{equation*}
Moreover, we have
\begin{equation*}
  \text{(III)} \quad  \sum_{k< b}w_1 \omega_k < \varepsilon\sum_{k\ge b} n_k \omega_k \quad \text{\ and \ \quad (IV)}\quad
  \bar w_1+\sum_{k> a} w_1 \omega_k < \varepsilon\sum_{k\le a} n_k \omega_k.
\end{equation*}
Inequality~(III) follows from
\begin{equation}
    \sum_{k< b}w_1 \omega_k < t w_1 \bar w_1=\frac{t {w_1}^2}{w(N)} < \frac{\varepsilon m^\circ}{w(N)} \le \varepsilon\sum_{k\ge b} n_k \omega_k
\end{equation}
using  $1/w(N)\le \omega_k\le \bar w_1$ and (\ref{ie_intermed_step}).
Similarly, 
(IV) follows from
\begin{equation}
    \bar w_1 +\sum_{k> a} w_1 \omega_k\le 
\bar w_1 +(t-1) w_1 \bar w_1 \le t \frac{{w_1}^2}{w(N)}<\varepsilon m^\circ\frac{1}{w(N)}\le \varepsilon\sum_{k\le a} n_k \omega_k.
\end{equation}
Let $n_k^\M$ denote the number of $k$-type voters in a coalition $\M\subseteq N$ and define
\begin{equation}
    D(\M)=\frac{\sum\limits_{k\le a} n_k^\M \omega_k}{\sum\limits_{k\le a} n_k \omega_k}\quad\text{and}\quad
    I(\M)=\frac{\sum\limits_{k\ge b} n_k^\M \omega_k}{\sum\limits_{k\ge b} n_k \omega_k}.
\end{equation}
$D(\M)$ is the share of the total weight of the $a$ ``most over-represented'' types (all having maximal nucleolus-to-relative weight ratio $x^*_1/\omega_1$) which they contribute in coalition $\M$.
Similarly, $I(\M)$ is the respective share for the $t-b+1$ ``most under-represented'' types.

Given a suitably large coalition $S\subseteq N$, replacing $w_h$ members of type~$k$~-- all with absolute weight $w_k$~-- by $w_k$ players of type $h$ yields a coalition $S'$ with $w(S')=w(S)$.
But if $r_k>r_h$, such replacement yields $x^*(S')<x^*(S)$.
Thus, for a MWC $\M$ with \emph{maximum} excess at $x^*$, i.e., with excess $1-x^*(T)\ge v(S)-x^*(S)$ for all $S\subseteq N$, there are no $k,h$ with $r_k>r_h$ such that (i) $w_1$ or more type~$k$-players belong to $\M$ and (ii) $w_1$ or more type~$h$-players do \emph{not} belong to $\M$.
This consideration restricts the numbers of members $n_k^\M$ of players of type~$k$ in any MWC $\M$ with maximum excess. There are three cases, for each of which we show $I(\M)-D(\M)> 
0$:
\begin{description}\item[\textnormal{Case 1:}]
$n_k^\M < w_1$ for all types $1\le k< b$. \smallskip

Then the relative weight $\sum_{k\le a} n_k^\M \omega_k$ in $\M$ of the most over-represented types is less than $\sum_{k\le a} w_1 \omega_k$. So inequality~(I) implies $D(\M)<\varepsilon$.
Since $\M$ is a winning coalition, the weight $\sum_{k\ge b} n_k^\M \omega_k$ in $\M$ of the most under-represented types is greater than $\bar q-\sum_{k<b} w_1 \omega_k$.
Due to (III) and $\sum_{k\ge b} n_k  \omega_k\le 1$ we have $I(\M)> \bar q-\varepsilon$.
So $I(\M)-D(\M)> \bar q-2\varepsilon \ge 0$.

\item[\textnormal{Case 2:}]
$n_k^\M\ge w_1$ for some $1\le k\le a$ but $n_h- n_h^{\M}<w_1$ for all $a< h\le t$.\footnote{If Case~1 does not apply, there is a smallest index $1\le k<b$ with $n_k^T\ge w_1$. Assume $k\le a$ first. Because $r_k>r_h$ for all $a<h<t$, the number $n_h-n_h^T$ of $h$-types outside coalition $T$ is less than $w_1$: otherwise the indicated replacement would yield a MWC $T'$ with $x^*(T')<x^*(T)$,  contradicting the maximum-excess property of $T.$ This is the description of Case~2. The remaining Case~3 involves $a<k<b$ where $r_k>r_h$ for all $b\le h\le t$. Then, analogously, $n_h-n_h^T<w_1$ must hold.} \smallskip

Using that $T$ is a MWC, the relative weight $\sum_{k\le a} n_k^\M \omega_k$ in $\M$ of the most over-represented types is less than $\bar q+\bar w_1-\sum_{k>a} \left(n_k-w_1\right) \omega_k$ in this case. So
inequality~(IV) and $\sum_{k\le a} n_k \omega_k\le 1$ imply $D(\M)< \bar q+\varepsilon$.
Moreover, the weight $\sum_{k\ge b} n_k^\M \omega_k$ in $\M$ of the most under-represented types is greater than $\sum_{k\ge b} \left(n_k-w_1\right) \omega_k$ and inequality~(II) implies $I(\M)>1-\varepsilon$. So $I(\M)-D(\M)>1-\bar q-2\varepsilon \ge 0$.

\item[\textnormal{Case 3:}]
$n_l^\M < w_1$ for all $1\le l\le a$ and $n_k^\M\ge w_1$ for some $a< k< b$ but
       $n_h- n_h^{\M}<w_1$ for all $b\le h\le t$. \smallskip

The relative weight $\sum_{l\le a} n_l^\M \omega_l$ in $\M$ of the most over-represented types is then less than $\sum_{l\le a}w_1 \omega_l$.
So inequality~(I) implies $D(\M)<\varepsilon$.
Similarly, the total weight of the players of types $b\le h\le t$ is greater than $\sum_{h\ge b}\left(n_h-w_1\right) \omega_h$.
Inequality~(II) then implies $I(\M)>1-\varepsilon$ and we have $I(\M)-D(\M)> 1-2\varepsilon>\bar q-2\varepsilon\ge 0$.
\end{description}

Recall that $x^*_k\ge w_k$ for all $1\le k\le a$. So for sufficiently small $\sigma>0$ \begin{equation}\label{eq:newimputation}
 x^{**}_k(\sigma)=
    \begin{cases}
    x^{*}_k-\sigma \omega_k & \text{if } 1\le k \le a, \\
    x^{*}_k & \text{if } a< k < b, \text{ and} \\
    x^{*}_k+\delta \sigma \omega_k & \text{if } b\le k \le t \\
    \end{cases}
\end{equation}
with $\delta=\sum_{k\le a} n_k \omega_k / \sum_{l\ge b} n_l \omega_l>0$ is an imputation.
$x^{**}_k(\sigma)$'s continuity implies existence of $\sigma>0$ so that no $S$ with $e(S,x^*)<E_1(x^*)$ has maximum excess at $x^{**}(\sigma)$. We fix such a value of $\sigma$ and write $x^{**}=x^{**}(\sigma)$.

It then suffices to consider coalitions $T'$ with maximum excess at $x^*$ in order to show the contradiction $E_1(x^{**})<E_1(x^{*})$.
Such $T'$ has to be winning, and for any MWC $T\subseteq T'$ it must be true that $e(x^*,T)=e(x^*,T')=E_1(x^*)$.
Since $T$ and $T'$ both are winning we have $e(x^{**},T')\le e(x^{**},T)$ and
\begin{equation}\label{eq:E_1alsMax}
E_1(x^{**})=\max\{e(x^{**},T)\colon T\text{ is }MWC \text{ and } e(x^*,T)=E_1(x^*)\}.
\end{equation}
Moreover, for every $T$ on the right-hand side of equation (\ref{eq:E_1alsMax})
\begin{equation}
  e(x^{**},T)-E_1(x^{*})=e(x^{**},T)-e(x^{*},T)=-\sigma\cdot \big(\underbrace{I(T)-D(T)}_{>0}\big)\cdot \sum_{k\le a} n_k \omega_k<0
\end{equation}
implies $e(x^{**},T)<E_1(x^{*})$, so that $E_1(x^{**})<E_1(x^{*})$. \hfill{$\square$}

\bigskip


\newcommand{\noopsort}[1]{}

\end{document}